\documentstyle[12pt]{article}
\newtheorem{theorem}{THEOREM}[section]
\newtheorem{definition}[theorem]{Definition}

\newtheorem{proposition}[theorem]{Proposition}
\def\wbar{\overline{w}}

\def\hbar{\overline{h}}

\def\Omegabar{\overline{\Omega}}

\def\CC{{\rm\kern.24em\vrule
width.02em height1.4ex
depth-.05ex\kern-.26em C}}
\def\QQ{{\rm\kern.24em\vrule width.02em
height1.4ex depth-.05ex\kern-.26em Q}}
\def\RR{{\rm I\kern-.2em R}}
\def\HH{{\rm I\kern-.2em H}}
\def\ZZ{{\rm\kern.26em\vrule width.02em
height0.5ex depth0ex\kern.04em\vrule width.02em
height1.47ex depth-1ex\kern-.34em Z}}
\def\Ibb#1{{\rm I\kern-.23em#1}}
\def\Ib#1{{\rm I\kern-.25em#1}}
\def\k#1{\kern#1em}
\def\vb#1{\vrule width.02em height1.4ex depth-.05ex}
\def\NN{\Ibb N}
\def\11{{\rm\k{.45}\vb0\k{-.142}1}}

\def\epf{\hskip.2in\vrule width.4pt height6.65pt
depth.15pt\vrule
width2.5pt height6.65pt depth-6.25pt\hskip-2.5pt\vrule
width2.5pt
height.25pt depth.15pt\vrule width.4pt
height6.65pt depth.15pt\ }

\def\proof{\noindent {\bf Proof. } }

\def \Omegabar{\overline \Omega}

\def \d{\partial}

\def\zbar{\overline{z}}
\def\dbar{\overline{\partial}}
\def \hbar{\overline{h}}

\def\wbar{\overline{w}}
\def\phibar{\overline{\phi}}
\def\\GBB{\cal B}

\def\1bar{\overline{1}}
\def\2bar{\overline{2}}

\def\dist{\hbox{dist}}
\def\Dbar{\overline{D}}
\def\lambdabar{\overline{\lambda}}
\def\lambdabar{\overline{\lambda}}
\begin{document}

 \title{ Solving  the Kerzman's problem on the 
 sup-norm  estimate for $\dbar$ on product domains}
\author{Song-Ying Li}

\maketitle





					
	{\narrower\smallskip\noindent \noindent {\bf Abstract.}  In this paper, the author solves the long term open problem of Kerzman on sup-norm estimate for Cauchy-Riemann equation on polydisc
	in $n$-dimensional complex space. The problem has been open since 1971. He also extends and solves the problem on a bounded product domain $\Omega^n$,  where  $\Omega$ either is 
	simply connected with $C^{1,\alpha}$ boundary or satisfies a uniform exterior ball condition with piecewise $C^1$ boundary.
	\par \smallskip}

\section{Introduction}

Let $\Omega$ be a bounded pseudoconvex domain in $\CC^n$. Let $f\in L^2_{(0,1)}(\Omega)$ be any $\dbar$-closed $(0,1)$-form
with coefficients $f_j\in L^2(\Omega)$. By H\"omander's theorem \cite{H1}, there is a unique $u\in L^2(\Omega)$ with $u\perp \hbox{Ker}(\dbar)$ such that $\dbar u=f$.
The regularity theory  for Cauchy-Riemann equations became  a very important research area in several complex variables
for many decades. In particular, sup-norm estimate for $\dbar$
is the most difficult one. 
When $\Omega$ is a smoothly bounded strictly pseudoconvex domain in $\CC^n$,
 in 1970,  Henkin \cite{Hen1},  Grauart and Lieb \cite{GL} constructed a formula solution for $\dbar u=f$  satisfying $\|u\|_{L^\infty }\le C_{\Omega}\|f\|_{L^\infty_{(0,1)}}$. 
 In 1971, Kerzman \cite{Ker} improved the above result in \cite{Hen1} and  \cite{GL}, 
 he  proved that $\|u\|_{C^\alpha (\Omega)}\le C_{\alpha, \Omega}\|f\|_{L^\infty_{(0,1)}}$
 for any $0<\alpha <1/2$. In 1971, Henkin and Romanov \cite{HR} proved the sharp estimate: $\|u\|_{C^{1/2}(\Omega)}\le C_\Omega \|f\|_{L^\infty_{(0,1)}}$.
 Recently, X. Gong \cite{Gon} generalized Henkin and Romanov's results.  He reduced  the assumption $\d\Omega\in C^\infty$ to $\d \Omega\in C^2$ and proved 
 that $\|u\|_{C^{\gamma +1/2}(\Omega)}
 \le \|f\|_{C^\gamma_{(0,1)}(\Omega)}$ for any $\gamma$ with that $\gamma+1/2$ is not an integer.
 In \cite{Ker}, when $\Omega=D^n$ is the unit polydisc in $\CC^n$,  Kerzman asked the following question: {\it  Does $\dbar u=f$
 have a solution satisfying $\|u\|_{C^\alpha }\le C_\alpha\|f\|_{L^\infty_{(0,1)}}$ for some $\alpha>0$ ?}  Let
 $f_j(\lambda)\in L^\infty(D) $ be  holomorphic in $D$ such that $u_0=\zbar_1 f_1(z_2)+\zbar_2 f_2(z_1) \not\in C(\overline{D}^2)$.  
 Let $f(z) =f_1(z_2) d\zbar_1+f_2(z_1) d\zbar_2$. Then $\dbar f=0$
 and $u_0\in L^\infty(D^2) \setminus C(\overline{D}^2)$ with $u_0\perp \hbox{Ker}(\dbar)$ solves $\dbar u=f$.
 Then the Kerzman's question can be refined by: {\it Does $\dbar u=f$ have a solution $u$ satisfying $\|u\|_{L^\infty}\le C\|f\|_{L^\infty_{(0,1)}} ?$} 
 The problem was studied by Henkin \cite{Hen}, he proved that if $f\in C^1_{(0,1)}(\Dbar^2)$ is $\dbar$-closed, then $\dbar u=f$ has a solution
 $u$ satisfying estimate $\|u\|_{L^\infty}\le C\|f\|_{L^\infty_{(0,1)}}$, where $C$ is a scalar constant. Notice that a $\dbar$-closed form
 $f \in L^\infty_{(0,1)}(D^n)$
 can not  be approximated by $\dbar$-closed forms in $C^1_{(0,1)}(\overline{D}^n)$ in $L^\infty(D^n)$-norm. Henkin's result only partially answered Kerzman's 
 question and left  the Kerzman's question remanning open.

In \cite{Lan}, Landucci was able to improve the solution $u$ of $\dbar u=f$ in \cite{Hen} to the canonical solution which is the solution $u_0\perp \hbox{Ker}(\dbar)$. Recently,
Chen and McNeal \cite{CM} introduced a new  space ${\cal B}^p_{(0,1)}(D^n)$ of $(0,1)$ over $D^n$ which is  smaller than $L^p_{(0,1)}(D^n)$
 and proved $L^p$-norm estimates for $f\in {\cal B}^p_{(0,1)}(D^n)$ for $1<p\le \infty$. Their result
generalized Henkin's result. For a simple example, they reduced Henkin's assumption: $f=f_1 d\zbar_1+f_2 d\zbar_2 \in C^1_{(0,1)}(\overline{D}^2)$ to $f\in L^\infty_{(0,1)}(D^2)$
satisfying ${\d f_1 \over \d \zbar_2}\in L^\infty(D^2)$.  Dong, Pan and
Zhang \cite{DPZ} proved a very clean and pretty theorem: If $\Omega$ is any bounded domain in $\CC$ with $C^2 $ boundary
and  $f\in C_{(0,1)}(\Omegabar^n)$ is $\dbar$-closed, then the canonical solution $u_0$ of $\dbar u=f$ satisfies $\|u_0\|_{L^\infty}\le C\|f\|_{L^\infty_{(0,1)}}$.
However, $C_{(0,1)}(\Omegabar^n) $ is strictly smaller than $L^\infty_{(0,1)}(\Omega^n)$, the Kerzman's question remains open (see \cite{Li}).   

Main purpose of the current paper is to give a complete solution of the Kerzman's long open problem on 
the unit polydisc in $\CC^n$. More general, we will prove that the canonical solution $u$ satisfying estimate $\|u\|_{\infty}\le C\|f\|_\infty$
 on  the product domains $\Omega^n$ for two classes of bounded domains $\Omega\subset \CC$.  The main theorem is stated as follows.

\begin{theorem} Let $\Omega$ be either a simply connected domain in $\CC$ with $C^{1,\alpha}$ boundary with some $\alpha>0$ or
a bounded domain with piecewise $C^1$ boundary satisfying a uniform exterior ball condition. Let
$f\in L^\infty_{(0,1)}(\Omega^n)$ be $\dbar$-closed. Then the canonical solution $u_0$ of  $\dbar u=f$ is constructed and satisfies
$$
\|u_0\|_{L^\infty(\Omega^n)}\le C\|f\|_{L^\infty_{(0,1)}(\Omega^n)}.\leqno(1.1)
$$
\end{theorem}

More informations for $\dbar$-estimates, one may find  from the following references as well as the references therein. For examples,  Chen and Shaw \cite{CS}, Fornaess and Sibony \cite{FS}, 
Krantz \cite{Kran, K76}, Range \cite{Ran}, Range and Siu \cite{RS72, RS73}, Shaw \cite{S91}
and Siu \cite{Siu}. For product domains, one may also see  \cite{CS},   \cite{DLT},  \cite{KL2} and other related articles in the reference.

The paper is organized as follows. In section 2, we provide  a formula solution for canonical solution of $\dbar u=f$  on the product domains.
In Section 3, technically, we translate the formula in Section 2 to one, from which we can get a uniform $L^p$ estimates. In Section 4,  we will prove Theorem 1.1. 
Finally, in Section 5, based on $\dbar$ -estimate  on the disc $D\subset \CC$, we give a  sharp theorem (Theorem 5.1) which is better  than Theorem 1.1.
\medskip

{\bf Acknowledgment.} The author would like to thank R-Y. Chen who read  through the first draft of manuscript and  Sun-sig Byun for providing some useful reference on Green's function.

\section{Formula Solutions}

\subsection{Green's functions}

Let $\Omega$ be a bounded domain in $\CC$ and
let $G(\lambda,\xi)$ be the Green's function for the Laplace operator ${\d^2 \over \d z \d \zbar}={1 \over 4} \Delta $ on 
$\Omega$. Then the Green's operator $G$ is defined by
$$
G[f](z)=\int_\Omega G(z, w) f(w) dA(w)\leqno(2.1)
$$
and  $G[f]$ satisfies
$$
{\d^2 G[f] \over \d \lambda \d\lambdabar}(\lambda) =f(\lambda). \leqno(2.2)
$$
Let $ A^2(\Omega)$ be the Bergman space over $\Omega$ which is the holomorphic subspace of $L^2(\Omega)$.
Let ${\cal P}:L^2(\Omega)\to A^2(\Omega)$ be the Bergman projection. Then
$$
(I-{\cal P})f (z)=-\int_\Omega {\d G(z, w) \over \d z \d \wbar} f(w) dA(w).\leqno(2.3)
$$
By Theorem 0.5 in Jerison and Kenig \cite{JK}, {\it if $\d \Omega$ is Lipschitz, there is a $p_1>4$ such that the Green's operator $G: W^{-1, p}(\Omega)\to W^{1, p}(\Omega)$
is bounded for $p_1'<p<p_1$.} (2.3) implies that if $\d \Omega$ is Lipschitz, then ${\cal P}: L^p(\Omega)\to A^p(\Omega)$ is bounded for $p_1'<p<p_1$.
One may find further information on regularity of Bergman projections  in \cite{LL}.

 We need some properties of the Green's function and estimations on the Green's function and its derivatives  based on the regularity
 of $\d \Omega$. We recall  a definition. {\it We say that  a bounded domain $\Omega \subset \RR^n$  satisfies  a uniform exterior ball (disc) condition if there is a positive number $r$ such that
 for any $z_0\in \d \Omega$, there is $z_0(r)\in \RR^n\setminus \Omegabar$ such that $\overline{B(z_0(r), r)} \cap \Omegabar=\{z_0\}$, where $B(x, r)$ is ball in $\RR^n$
  centered at $x$ with radius $r$.}  It is easy to see that if $\d \Omega$ is $C^2$, then
 $\Omega$ satisfies a uniform  exterior (and interior) ball condition.
 
The following theorem on the Green's function  was proved by  Gr\"uter and Widman  \cite{GW} (Theorem 3.3)
which was also stated  as Theorem 4.5  in \cite{MM}.  

\begin{theorem} If $\Omega$ is a bounded domain  in $\RR^n$ which satisfies a uniform exterior ball condition, then its associated Green function  satisfies the
following five properties for all $x, y\in \Omega$:

(i) $|G(x, y)|\le  C d_\Omega( x) |x-y|^{1-n}$;

(ii) $|G(x,y)|\le C d_\Omega (x)  d_\Omega (y)  |x-y|^{-n}$;

(iii) $ |\nabla_x G(x,y)|\le C|x-y|^{1-n}$;

(iv) $|\nabla_x G(x, y)| \le C d_\Omega (y) |x-y|^{-n}$;

(v) $|\nabla_x \nabla_y G(x, y)|\le C|x-y|^{-n} $.

\noindent Here $C$ is a constant depending only on $\Omega$ and $d_\Omega(x)$ is distance from $x$ to $\d \Omega$.
\end{theorem}


Notice that  $\Omega $ having $ C^{1,\alpha}$ boundary with $\alpha \in (0,1)$ may not  satisfy  a uniform exterior ball condition.
 We will give a formula for the Green's function on a bounded simply connected domain  in $\CC$ with $C^{1,\alpha}$ boundary. 

Applying the argument by Kerzman \cite{Ker2} and regularity theorem (Theorem 8.34 in \cite{GT}), one can prove the following result. 

\begin{proposition} Let $\Omega$ be a bounded domains in $\CC$ with
$C^{1,\alpha}$ boundary for some $0<\alpha<1$. 

(i)  If $\psi: \Omega\to D(0,1)$ is a proper
holomorphic map, then $\psi\in C^{1,\alpha}(\Omegabar_1)$;

(ii) If $\phi: \Omega\to D(0,1)$ is biholomorphic, then the Green's function $G_\Omega$  for ${\d^2 \over \d z\d \zbar}$ in $\Omega$ is given by
$$
G_{\Omega}(z,w)={1\over \pi} \log \Big|{\phi(z)-\phi(w)\over 1-\phi(z)\phibar(w)}\Big|^2 \leqno(2.4)
$$
which satisfies (i)--(v) in Theorem 2.1.
\end{proposition}

\proof By Theorem 8.34 in \cite{GT}, if $g\in L^\infty(D)$,  then 
$$
\Delta u=g \hbox{ in } D, \quad u=0 \hbox{ on }\d D
$$
has a unique solution $u\in C^{1,\alpha}(\overline{D} )$. Let $g\in C^\infty_0 (D)$ be  a non-negative function on $D$ such that
 $\{z\in D: g(z)>0\}$
is a non-empty, relatively compact subset in $D$. Let $v(z)=u(\psi(z))$ be a function on $\Omega$ which solves the Dirichlet boundary problem:
$$
\cases{\Delta v(z)=g(\psi(z))|\psi'(z)|^2,  \ z\in \Omega, \cr
\quad v(z)=0, \qquad\qquad \quad\   z\in \d \Omega.\cr}
$$
By the elliptic theory (Theorem 3.34 in \cite{GT}), one has $v\in C^{1,\alpha}(\Omegabar)$. Then
$$
{\d v\over \d z}(z)={\d u\over \d w}(\psi(z))\psi'(z).
$$
Since $D$ satisfies an  interior ball condition, by Hopf's lemma, one has ${\d u\over \d w}(w)\ne 0$ on $\d D$. Since
$u\in C^{1,\alpha}(\overline{D})$, one has ${\d u\over \d w}(w)\ne 0$ on the closed annulus
$A(0, 1-\epsilon, 1]=
\{w\in D: 1-\epsilon \le |w| \le 1\}$
for some small $\epsilon>0$. This implies
$$
\psi'(z) ={\d v(z) \over \d z} /{\d u\over \d w}(\psi(z))\quad \hbox{ on } \psi^{-1}(A(0, 1-\epsilon, 1]). \leqno(2.5)
$$
 This implies that $\psi \in C^1(\Omegabar)$ since $\psi$ is holomorphic in $\Omega$. Applying (2.5) again, one can see that 
$\psi'(z)\in C^\alpha (\Omegabar)$. Therefore, $\psi \in C^{1,\alpha}(\Omegabar)$.

\medskip
It is well known that the Green's function for ${\d^2\over \d z \d \zbar}$ in the unit disc $D$ is:
$$
G(z,w)={1\over \pi} \log\Big |{w-z\over 1-z\wbar}\Big |^2,\quad z, w\in D.\leqno(2.6)
$$
If $\phi: \Omega\to D$ is a bilomorphic map, then it is easy to check that the Greens's function for $\Omega$ is given by (2.4).
Moreover, one can check
that $G_\Omega$ satisfies Properties (i)--(v) in Theorem 2.1 when $n=2$. \epf

\subsection{Formula solution to $\dbar$-equations}
Let $G=G_\Omega$ be the Green's function for ${\d^2 \over \d z\d \zbar}$ on $\Omega$. Define
$$
k(z, w)={\d G_{\Omega}(z, w)\over \d z}
 \leqno(2.7)
$$
and
$$ \label{S1}
T[f](z)=\int_\Omega k(z,w) f(w) dA(w).\leqno(2.8)
$$

For simplicity, we give the following definition.
\begin{definition} A domain $\Omega\subset \CC$ is said to be admissible if either $\Omega$ is bounded, simply connected with $C^{1,\alpha}$ boundary
for some $\alpha\in (0,1)$ or $\Omega$ is bounded with piecewise $C^1$ boundary and satisfies a uniform exterior ball condition.
\end{definition}

\begin{proposition} Let $\Omega \subset \CC$ be an admissible domain and $2<p<\infty$. Then

(i)  If $f\in L^2(\Omega)$, then $T[f]  $ is 
 the canonical solution of  $\dbar u=f d \zbar$;
 
 (ii)  $T: L^p(\Omega)\to L^\infty(\Omega)$ is bounded;
 
 (iii) $T: L^p(\Omega)\to C^{1-2/p}(K)$  for any compact set $K\subset 
 \Omega$;
 
 (iv) If $\Omega$ is simply connected and $\d \Omega\in C^{1,\alpha_0}$, then 
 $T: L^p(\Omega)\to C^{\alpha}(\Omegabar)$, where $\alpha=\min\{\alpha_0, 1-2/p\}$.
\end{proposition}

\proof By (2.6) and (2.7), the definition of $T[f]$ and the definition of the Green's function, one can easily see that 
$$
{\d T [f]\over \d \lambdabar} (\lambda)={\d^2 G[f] \over \d \lambda\d \lambdabar} =f(\lambda),\quad \lambda \in \Omega.
$$
For any  $h(\lambda)\in W^{1,2}(\Omega)\cap A^2(\Omega)$ and Theorem  2.1, one has
\begin{eqnarray*}
\int_\Omega T[f](\lambda) \hbar (\lambda)dA(\lambda)&=&\int_\Omega \int_\Omega k(\lambda, w) \hbar(\lambda) dA(\lambda) f(w) dA(w)\\
&=-&\int_\Omega \int_\Omega G(\lambda, w) {\d \hbar(\lambda) \over \d \lambda} dA(\lambda) f(w) d A(w)\\
&=&-\int_\Omega 0\cdot  f(w) dA(w)\\
&=& 0.
\end{eqnarray*}
Since $W^{1,2}(\Omega)\cap A^2(\Omega)$ is dense in $A^2(\Omega)$, one has proved that
 $T[f]\perp A^2(\Omega)$. So, $T[f]$ is the canonical solution of $\dbar u =f d\zbar $ in $\Omega$. Part (i) is proved.
 
 For Part (ii),  by Part (iv) in Theorem 2.1,   Proposition 2.2 and (2.4),   one has
 $$
 |k(z, w)|=|{\d G(z, w) \over \d z}| \le {C\over |z-w|}.\leqno(2.9)
 $$
 This implies
 $$
 |T[f](z)|\le C\int_\Omega{ |f(w)|\over |w-z|} dA(w)\le {C\over 2-p'} \|f\|_{L^p} \le C{p-1\over p-2}  \|f\|_{L^p},
 $$
for any $2<p\le \infty$. This means $\|T[f]\|_{L^\infty}\le C{p-1\over p-2} \|f\|_{L^p}$ if $p>2$. Part (ii) is proved.
Let
$$
v(z)={1\over \pi}\int_\Omega {f(w)\over z-w} dA(w).
$$
Then ${\d v\over \d \zbar}=f$. By Sobolev embedding theorem, one has that  $v\in W^{1, p}(\Omega)\subset C^{1-2/p}(\Omegabar)$ for $2< p<\infty$.  Thus,
$$
T[f]=v-{\cal P}[v]\in C^{1-2/p}(K),\quad\hbox{for any compact set } K\subset \Omega.
$$
Therefore, Part (iii) is completed.

When $\Omega$ is simply connected and if $\phi: \Omega\to D(0,1)$ is a biholomorphic map, then Bergman kernel for $\Omega$ is
$$
K(z,w)={1\over \pi} {\phi'(z) \overline{\phi'(w)} \over (1-\phi (z) \phibar(w) )^2 },\quad z, w  \in \Omega.  \leqno(2.10)
$$
It is easy to verify that ${\cal P}[v] \subset C^\alpha(\Omegabar)$ with $ \alpha=\min\{\alpha_0, 1-2/p\}.$ This proves Part (iv). 
Therefore, the proof of the proposition is complete.\epf

\medskip

For any  $1\le j\le n$ and $z\in \CC^n$, write
$$
z^{(j)}=(z_1,\cdots, z_{j-1}, z_{j+1}, \cdots, z_n), \quad z=(z_j; z^{(j)}).\leqno(2.11)
$$
Let $f\in L^2(\Omega^n)$, we define
the Bergman projection $P_j: L^2(\Omega)\to A^2(\Omega)$  by
$$
P_j f (z)={\cal P} [f(\cdot, z^{(j)})](z_j)
=\int_\Omega K(z_j, w_j) f(w_j; z^{(j)}) dA(w_j), \leqno(2.12)
$$
for almost every $z^{(j)}\in \Omega^{n-1}.$ We also use the notations $P_0=P_{n+1}=I$. Similarly,  we also use the following notation:
$$
T_j f(z)=T[f(\cdot\, ; z^{(j)})] (z_j), \quad 1\le j\le n.\leqno(2.13)
$$

The following theorem is a very important formulation for the canonical solution of $\dbar u=f$.

\begin{theorem} Let $\Omega$ be an admissible domain in $\CC$. For $2< p\le \infty$ and any $\dbar$-closed $(0,1)$-form
 $f=\sum_{j=1}^n f_j d\zbar_j \in L^p_{(0,1)}(\Omega^n)$, the canonical solution 
$u=S[f] \in L^2(\Omega^n)$  to $\dbar u=f$ satisfies
$$
S[f](z)=\sum_{j=1}^n T_j  P_{j-1}\cdots P_{0} f_j=\sum_{j=1}^n T_j P_{j+1}\cdots P_{n+1} f_j. \leqno(2.14)
$$
\end{theorem}

\proof  For each $1\le j\le n$, since ${\d u(z_j; z^{(j)})\over \d \zbar_j} =f_j ( z_j; z^{(j)})\in L^p(\Omega)$.
By the estimates on the Green's function given  by Theorem 2.1, Propositions 2.2 and 2.4, one has that
$$
u(z_j; z^{(j)})-P_j[u(\cdot\, ; z^{(j)})](z_j)=T_j[f_j(\cdot\, ; z^{(j)})](z_j),\leqno(2.15)
$$
for almost every $ z^{(j)}\in \Omega^{n-1}.$ 

  Since $u-P_1[u]$ is the canonical solution of ${\d u\over \d \zbar_1}=f_1$, one has
$$
P_0 u-P_1 P_0 u=u-P_1[u]=T_1 f_1=T_1 P_0 f_1.
$$
Similarly, $P_1P_0[u]-P_2P_1P_0[u]=T_2 P_1 f_2$. Keeping the same process, one has
$$
P_{j-1}\cdots P_1 P_0u-P_j P_{j-1} \cdots P_1P_0 u= T_j P_{j-1}\cdots P_1  f_j,\quad 1\le j\le n.
$$
Since $P_1\cdots P_n u=0$ and $P_0=I$, one has
$$
S[f]=u=\sum_{j=1}^n (P_{j-1}\cdots P_0u-P_j P_{j-1} \cdots P_0 u )=\sum_{j=1}^n  T_j P_{j-1}\cdots P_{1} P_0 f_j.
$$

On the other hands, let $P_{n+1}=I$, then
$$
 u-P_n  u=T_n f_n.
$$
With the same process, one has
$$
P_n \cdots P_j u-P_n\cdots P_j  P_{j-1}  u=T_{j-1} P_j\cdots P_n f_{j-1}.
$$
Since $u$ is the canonical solution of $\dbar u=f$, one has $P_{n+1} P_n\cdots P_1 u=0$ and
$$
\sum_{j=1}^{n} T_{j} P_{j+1}\cdots P_{n+1} f_{j}=\sum_{j=1}^{n} (P_{n+1} P_n \cdots P_{j+1} u- P_{n+1} P_n\cdots P_j   u)=u.
$$
These prove (2.14), so, the proof of  Theorem 2.5  is complete. \epf

If $\Omega$ is a simply connected domain with $C^{1,\alpha}$ boundary. Let $\phi: \Omega\to D$ be a biholomorphic mapping. Then 
the Bergman kernel function is given by (2.10).
Since $\phi\in C^{1,\alpha}(\Omegabar)$, one has that the Bergman projection $P: L^p(\Omega)\to L^p(\Omega)$ is bounded for all
$1<p<\infty$. By the expression of $S[f]$, one can easily see the following statement holds.

\begin{theorem}  Let $1<p<\infty$ and  let $\Omega$ be a bounded simply connected domain in $\CC$ with $C^{1,\alpha}$ boundary for some $\alpha>0$.
View $S[f]$ as a linear operator on $L^p_{(0,1)}(\Omega^n)$ defined by (2.14).
If   $f_m , f\in L^p_{(0,1)}(\Omega^n)$ with $f_m\to f$ in $L^p_{(0,1)}(\Omega)$, then
$$
\lim_{m\to \infty}\Big \|S[f_m]-S[f]\Big\|_{L^p_{(0,1)}(\Omega^n)}=0.\leqno(2.16)
$$
\end{theorem}

When $\Omega$ is a bounded domain with piecewise $C^1$ boundary and satisfies  a uniform exterior ball condition, we don't know whether the Bergman
projection $P:L^p(\Omega)\to A^p(\Omega)$ is bounded or not for all $4<p_1\le p<\infty$. However,  with the different expression of $S[f]$ 
 given in the next section,
we will be able to prove Theorem 2.6 remains true under the assumtion $\dbar f_m=0$ and $\dbar f=0$.

\begin{theorem}  Let $1<p<\infty$ and  let $\Omega$ be a bounded domain in $\CC$ with piecewise $C^1$ boundary satisfying a uniform exterior ball condition.
If   $f_m \in C^1_{(0,1)}(\Omegabar^n)$ and $f\in L^p_{(0,1)}(\Omega)$ are $\dbar$-closed and satisfy $f_m\to f$ in $L^p_{(0,1)}(\Omega^n)$ as $m\to \infty$, then
$$
\lim_{m\to \infty}\Big \|S[f_m]-S[f]\Big\|_{L^p_{(0,1)} (\Omega^n)}=0
\quad\hbox{and}\quad
\dbar S[f]=f. \leqno(2.17)
$$
\end{theorem}


\section{Regularity and  a new formula solution}
For any $1\le i \ne j\le n$, define
$$
\tau_{i,j}(z,w)=|w_i-z_i|^2+ | w_j-z_j|^2=\tau_{j,i}(z,w)\leqno(3.1)
$$
and
\begin{eqnarray*}
(3.2)\qquad b^{i, j}(z,w)&:=&{\d \over \d \wbar_j} ({| w_j-z_j|^2 k(z_j, w_j)\over \tau_{i,j}(z,w)})\\
&=& k(z_j, w_j) {\d\over \d \wbar_j} ({|w_j-z_j|^2\over \tau_{i,j}})+{|w_j-z_j|^2\over \tau_{i,j}}{\d k(z_j, w_j)\over \d \wbar_j} \qquad\qquad\\
&=& k(z_j, w_j) {(w_j-z_j) |w_i-z_i|^2\over \tau_{i,j}^2 }+{|w_j-z_j|^2\over \tau_{i,j}}{\d k(z_j, w_j)\over \d \wbar_j}\\
&=&h(z_j, w_j){|w_i-z_i|^2\over \tau_{i,j}^2}+{H(z_j, w_j)\over \tau_{i,j}},
\end{eqnarray*}
where
$$
h(z_j, w_j)=(w_j-z_j)k(z_j, w_j), \hbox{ and } H(z_j, w_j)=|w_j-z_j|^2{\d k(z_j, w_j) \over \d \wbar_j}.\leqno(3.3)
$$
By Theorem 2.1 and Proposition 2.2, with $C=C_\Omega$, one has
$$
|h(z_j, w_j|+ |H(z_j, w_j)|\le C  \ \hbox{ and } \  |h(z_j, w_j)|\le {C d_\Omega(w_j) \over |z_j-w_j|}. \leqno(3.4)
$$
Therefore
$$
|b^{i,j}(z, w)|\le {C\over \tau_{i,j}(z,w)}.\leqno(3.5)
$$
Notice that
$$
{\d b^{j,i}\over \d \wbar_j}
=h(z_i, w_i) {(w_j-z_j)(|w_i-z_i|^2-|w_j-z_j|^2)\over \tau_{ij}^3}-H(z_i, w_i){w_j-z_j\over \tau_{i,j}^2}.\leqno(3.6)
$$
Then
$$
\Big|{\d b^{j,i}\over \d \wbar_j}\Big|\le C{ |w_j-z_j| \over \tau_{i,j}^2}.\leqno(3.7)
$$
Write
\begin{eqnarray*}
\lefteqn{(3.8)\qquad\qquad \qquad{\d \over \d \wbar_j} \Big( b^{j, i}(z,w) {|w_j-z_j|^2\over \tau_{j,k}(z,w)} k(z_j, w_j))}\\
&\qquad\qquad \qquad\qquad =&b^{j, i} b^{k, j}
+{|w_j-z_j|^2 \over \tau_{j, k}}k(z_j, w_j){\d b^{j,i}\over \d \wbar_j} \qquad\qquad \qquad\qquad\quad\\
&\qquad\qquad \qquad\qquad  =&b^{j,i} b^{k, j} +{a^{j, i}\over \tau_{j,k}},
\end{eqnarray*}
where
$$
a^{j, i}=|w_j-z_j|^2 k(z_j, w_j){\d b^{j,i}\over \d \wbar_j},\quad |a^{j, i}|\le C {|w_j-z_j|^2\over \tau_{i,j}^2}\le {C\over \tau_{i,j}}. \leqno(3.9)
$$
Let
$$
B_{j , i} [g]=\int_{\Omega} g(w) b^{j, i}(z, w) dA(w_i)\leqno(3.10)
$$
and
$$
A^k_{ j, i}[g]=\int_{\Omega^2} { a^{j, i}\over \tau_{j, k} } g(w) dA(w_i) dA(w_j).\leqno(3.11)
$$
\begin{proposition} Let $f\in C^1_{(0,1)}(\Omegabar)$ be $\dbar$-closed.  Then
for any $i\ne j$, one has
$$
T_jT_i[{\d f_j \over\d \zbar_i}]=-T_jB_{j,i}[f_j]-T_i B_{i,j}[f_i], \leqno(3.12)
$$
$$
T_jP_i[f_j]=T_j[f_j] -T_jT_i[{\d f_j\over \d \zbar_i}]=T_j[f_j]+T_j B_{j,i}[f_j]+T_i B_{i,j}[f_i]\leqno(3.13)
$$
and
$$
T_i T_j B_{j,k}[{\d f_j\over \d \zbar_i} ]=-T_j B_{j ,i} B_{j,k}[f_j] -T_i B_{i, j} B_{j, k}[f_i]-T_i  A^i_{j, k}[f_i] .\leqno(3.14)
$$
\end{proposition}
\proof  Since $f$ is $\dbar$-closed, one has
$$
{\d f_j \over \d \zbar_i}={|w_i-z_i|^2\over \tau_{j,i}(z,w)}{\d f_j \over \d \zbar_i}+{|w_j-z_j|^2\over \tau_{i,j}(z,w)}{\d f_i \over \d \zbar_j}. \leqno(3.15)
$$ 
Notice that  $|k(z_i, w_i)| |w_i-z_i|^2\le C d_\Omega(w_i)$ and integration by part, one has
\begin{eqnarray*}
T_j T_i[{\d f_j\over \d \zbar_i}]
&=&\int_{\Omega^2}  k(z_i, w_i) k(z_j, w_j) {\d f_j\over \d \wbar_i} {|w_i-z_i|^2\over \tau_{i,j}(z,w)}dA(w_i) dA(w_j)\\
&&+\int_{\Omega^2}  k(z_i, w_i) k(z_j, w_j) {\d f_i\over \d \wbar_j } {|w_j-z_j|^2\over \tau_{i,j}(z,w)}dA(w_j) dA(w_i)\\
&=&-T_j B_{j, i} [f_j] -T_i B_{i,j}[f_i].
\end{eqnarray*}
(3.12) is proved. Since
$$
T_j P_i[f_j]=T_j[f_j]-T_j(I-P_i) f_j=T_j [f_j]- T_j T_i\Big[{\d f_j\over \d \zbar_i}\Big],
$$
by (3.12), one has proved  (3.13).  For simplicity, if no confusions may cause,  we let
$$
k_j=k(z_j, w_j),\quad 1\le j\le n.
$$
Then
$$
k_i k_j b^{j, k} {\d f_j\over \d \wbar_i}
=k_j  b^{j, k} {\d f_j\over \d \wbar_i} k_i {|w_i-z_i|^2\over \tau_{i,j}}+k_i b^{j, k} {\d f_i\over \d \wbar_j}{|w_j-z_j|^2\over \tau_{i,j}} k_j.
$$
By (3.8), 
$$
{\d\over \d \wbar_j} [b^{j, k}  {|w_j-z_j|^2\over \tau_{i,j}} k_j]=b^{j, k} b^{i, j}+{a^{j, k}\over \tau_{i,j}}.\leqno(3.16)
$$
By (3.8)--(3.11) and integration by part, one has
$$
 -T_iT_j B_{j,k}[{\d f_j\over \d \zbar_i}]
= T_j B_{j,i} B_{j,k} [f_j] +T_i B_{i,j} B_{j, k}[f_i]+T_i A^i_{j, k}[f_i] .\leqno(3.17)
$$
Therefore,  (3.14) is proved, so is the proposition. \epf
\medskip

Write 
$$
I=(i_1,i_1,\cdots, i_k) \hbox{ with } 1\le  i_1<i_2<\cdots< i_k\le n.  \leqno(3.18)
$$
For each $1\le \ell \le n$, we let $I=(i_1,\cdots, i_k) $ with $i_j\in \{1,\cdots, n\}\setminus\{\ell \}$ for $1\le j\le k$. Let
$E^\ell_I (z,w)$ be  an integrable function in $(z_\ell , z_{i_1},\cdots, z_{i_k})$ and  in $(w_\ell, w_{i_1},\cdots, w_{i_k})$ over $\Omega^{k+1}$
 satisfying the estimate: 
$$
 |E^\ell _I(z,w)|\le {C\over |w_\ell -z_\ell|^{1+k\epsilon } \ell_I(\epsilon)},\quad
\ell_I(\epsilon)=:\prod_{j=1}^k  |w_{i_j}-z_{i_j}|^{2-\epsilon}\leqno(3.19)
$$
for any small $\epsilon>0$. 

For each $I \subset \{1,\cdots, n\}\setminus \{\ell \}$ with $|I|=k$, we define
$$
T^\ell _I[f_i] =\int_{\Omega ^{k+1}}  E^\ell_I (z, w) f_i(w) d v (w_\ell, w_{i_1},\cdots, w_{i_k}).\leqno(3.20)
$$
We are going to prove the following theorem.

\begin{theorem} Let $f\in C^1_{(0,1)}(\Omegabar)$ be $\dbar$-closed. Then there exist $E^j_I$ satisfy (3.19) and $T^j_I$ defined by (3.20) such that
$$
S[f](z)=\sum_{j=1}^n T_j[f_j]+\sum_{j=1}^n  \sum_{|I| \le n-1} T^j_I[f_j]. \leqno(3.21)
$$
\end{theorem}

\proof It is  obvious if $n=1$. We start with $n=2$. Since (2.12) and (3.13), one has
$$
S[f]=T_1[f_1]+T_2P_1[ f_2]=T_1 f_1+T_2[f_2]+T_2 B_{2,1}[f_2]+T_1B_{1,2}[f_1].
$$
Then 
$$
E^2_1=k(z_2, w_2) b^{2,1} \quad\hbox{and } \quad E^1_2= k(z_1, w_1) b^{1,2}.
$$
Applying 
$$
a^{\epsilon} b^{2-\epsilon} \le {\epsilon \over 2} a^2+ {2-\epsilon \over 2} b^2\le a^2+b^2 \leqno(3.22)
$$
and estimate (3.5) on $b^{i, j}$, one has
$$
|E^2_1(z, w)|\le {C\over |w_2-z_2|} {C\over \tau_{1,2}}\le {C\over |w_2-z_2|^{1+\epsilon} |w_1-z_1|^{2-\epsilon}}.
$$
Similarly,
$$
|E^1_2(z, w)|\le {C\over |w_1-z_1|^{1+\epsilon} |w_2-z_2|^{2-\epsilon}}.
$$
This prove the case $n=2$. 

For any $i< j<k$, notice that $(I-P_j)[f_k]=T_j[ {\d f_k \over \d\zbar_j}]$,  one has 
\begin{eqnarray*}
(3.23)\qquad T_k P_j P_i [f_k]&=&T_k P_i [f_k]-T_k T_j P_i [{\d f_k\over \d \zbar_j}]\\
&=&T_k[ f_k] -T_k T_i[{\d f_k \over \d \zbar_i}]
-P_i T_k T_j  [{\d f_k \over \d \zbar_j}]\qquad\qquad\qquad\quad
\end{eqnarray*}
and
\begin{eqnarray*}
\lefteqn{ -P_i T_k T_j  [{\d f_k \over \d \zbar_j}]}\\
&\qquad =&P_i T_k B_{k,j}[f_k]+P_i T_j B_{j, k}[f_j]\\
&\qquad =&
T_k B_{k,j}[f_k]+T_j B_{j, k}[f_j]-T_i T_k B_{k,j}[{\d f_k\over \d \zbar_i} ]- T_iT_j B_{j, k}[{\d f_j\over \d \zbar_i} ].\qquad
\end{eqnarray*}
Therefore, combining (3.12),  (3.14),  (3.20),  (3.21) and the above, one has
\begin{eqnarray*}
(3.24)\quad T_k P_j P_i[f_k]&=&T_k [f_k]+T_k B_{k,i}[f_k]+T_i B_{i, k}[f_i]+
T_k B_{k,j}[f_k]+T_j B_{j, k}[f_j]\\
&+& T_j B_{j, i} B_{j, k}[f_j]+T_i B_{i, j} B_{j, k}[f_i] +T_i A^i_{ j, k}[f_i]\\
&+&T_k B_{k,i} B_{k,j}[f_k] +T_i B_{i, k}  B_{k, j}[f_i ]+T_i  A^i_{k, j}[f_i] .
\end{eqnarray*}
By (3.5) and (3.22), one has
$$
|E^j_{i, k}|=|k(z_j, w_j)b^{j, i} b^{j,k}|
\le {C\over |w_j-z_j| \tau_{j, i} \tau_{j, k}}\le {C\over |w_j-z_j|^{1+2\epsilon} \ell_{i, j}(\epsilon)}. \leqno(3.25)
$$
Similarly,
$$
|E^k_{i,j}|\le {C\over |w_k-z_k|^{1+2\epsilon} \ell_{i,j}(\epsilon)}.\leqno(3.26)
$$
By (3.5), (3.9) and (3.22), one has
\begin{eqnarray*}
(3.27)\quad |E^i_{j, k}|&=& |k(z_i, w_i)| \Big |[b^{i,j} b^{j, k}+b^{i, k} b^{k, j}+{a^{j,k} \over \tau_{i, j}}
+{a^{k, j}\over \tau_{i, k}}]\Big |\\
&=&{C\over |w_i-z_i| } \Big({1 \over \tau_{i, j} \tau_{j,k}}+{1\over \tau_{i, k} \tau_{k, j}}
+{1 \over \tau_{j, k} \tau_{i, j} }+{1 \over \tau_{k,j}\tau_{i,k}}\Big)\\
&\le&{C\over |w_i-z_i|^{1+2\epsilon} \ell_{j,k}(\epsilon) }.
\end{eqnarray*}
By (3.24)--(3.27), (3.19) and Theorem 2.5, we have proved Theorem 3.2 when $n=3$.

Notice that for $k\ge 4$, one has
$$
T_kP_{k-1}\cdots P_1[f_k]=T_k P_{k-1}\cdots P_2 [f_k]-P_2\cdots P_{k-1} T_k T_1 [{\d f_k\over \d \zbar_1}]\leqno(3.28)
$$
and by (3.12)
$$
(3.29)\quad
-P_2\cdots P_{k-1}T_k T_1 [{\d f_k\over \d \zbar_1}]
=P_2\cdots P_{k-1} T_k B_{k,1} [f_k] +P_2\cdots P_{k-1} T_1 B_{1, k}  [f_1].
$$
One may use the principle of mathematics induction to complete the proof of Theorem 3.2. We continue to demonstrate the case $k=4$. By (3.24)
and (3.28)--(3.29), one need only to consider $P_2\cdots P_{k-1} T_k B_{k,1} [f_k]$, the other term in (3.29) can be computed similarly by exchange $k$ and $1$.
By (3.13), one has
$$
 P_2P_3 T_k B_{k,1} [f_k]
=T_k B_{k, 1}[f_k] -T_3 T_k B_{k, 1}[{\d f_k \over \d \zbar_3}]-P_3 T_2 T_k B_{k, 1}[{\d f_k\over \d \zbar_2}]. \leqno(3.30)
$$
By (3.14), one has
$$
T_i T_k B_{k,1}[{\d f_k\over \d \zbar_i} ]=-T_k B_{k ,i} B_{k,1}[f_k] -T_i B_{i, k} B_{k, 1}[f_i]-T_i  A^i_{k, 1}[f_i] \leqno(3.32)
$$
and
\begin{eqnarray*}
\lefteqn{(3.33)\qquad -P_3 T_2 T_k B_{k, 1}[{\d f_k\over \d \zbar_2}]}\\
&\qquad\qquad =&P_3 T_k B_{k ,2} B_{k,1}[f_k] +P_3 T_2 B_{2, k} B_{k, 1}[f_2]+P_3T_2 A^2_{k, 1}[f_2] \\
&\qquad\qquad =&T_k B_{k ,2} B_{k,1}[f_k] +T_2  B_{2, k} B_{k, 1}[f_2] +T_2 A^2_{k, 1}[f_2] \\
& \qquad \qquad& -T_k T_3 B_{k ,2} B_{k,1}[{\d f_k \over \d \zbar_3}] - T_2 T_3B_{2, k} B_{k, 1}[{\d f_2 \over \d \zbar_3}]- T_2 T_3 A^2_{k, 1}[{\d f_2\over \d \zbar_3}].\qquad\qquad 
\end{eqnarray*}
By (3.16), one has
$$
{\d\over \d \wbar_k} [b^{k, 1}  {|w_k-z_k|^2\over \tau_{3 ,k}} k(z_k, w_k)]=b^{k, 1} b^{3, k}+{a^{k, 1}\over \tau_{3,k}}\leqno(3.34)
$$
and
$$
{\d\over \d \wbar_k} [b^{k, 2}  {|w_k-z_k|^2\over \tau_{3,k}} k(z_k, w_k)]=b^{k, 2} b^{3, k}+{a^{k, 2}\over \tau_{3,k}}.\leqno(3.34')
$$
Then
\begin{eqnarray*}
\lefteqn{ (3.35)\  -T_k T_3 B_{k ,2} B_{k,1}[{\d f_k \over \d \zbar_3}] }\\
& \qquad =&T_k B_{k,2} B_{k,1} B_{k, 3} [f_k]+T_3 B_{k, 2} B_{k, 1} B_{3, k} [f_3]
+T_3 B_{k,2} A^3_{k,1} [f_3]\\
&&
 +T_3 B_{k, 1} B_{k, 2} B_{3, k} [f_3]+T_3 B_{k, 1} A^3_{k,2} [f_3]\\
 & \qquad =&T_k B_{k,2} B_{k,1} B_{k, 3} [f_k]+2T_3 B_{k, 2} B_{k, 1} B_{3, k} [f_3]
+T_3 B_{k,2} A^3_{k,1} [f_3]
+T_3 B_{k, 1} A^3_{k,2} [f_3].
\end{eqnarray*}
Since
$$
 {\d \over \d\wbar_2} {a^{k, 1} \over \tau_{2, k}}
= -{a^{k, 1} (w_2-z_2) \over \tau_{2, k}^2},
$$
one has
\begin{eqnarray*}
\lefteqn{(3.36)\  -T_2 T_3 A^2_{k,1} [{\d f_2\over \d \zbar_3}]}\\
& \qquad=&T_2 A^2_{k, 1} B_{2,3} [f_2]+ T_3 A^2_{k,1}B_{3, 2}[f_3]\\
&&+T_3[ \int_{\Omega^3} k(z_2, w_2) {|w_2-z_2|^2\over \tau_{2,3}} {a^{k, 1} (z_2-w_2) \over \tau_{2, k}^2} f_3 dA(w_1) dA(w_2) dA(w_k) ].
\end{eqnarray*}
Write
$$
a^{2, k, 1}(z,w)=k(z_2, w_2) |w_2-z_2|^2  {a^{k, 1} (z_2-w_2) \over \tau_{k, 2}^2} \leqno(3.37)
$$
and
$$
A^3_{2, k,1} [f_3]=\int_{\Omega^3} {a^{2, k,1} (z,w) \over \tau_{2, 3}} dA(w_2) dA(w_k) dA(w_1).\leqno(3.38)
$$
Then
$$
|a^{2,k, 1}(z,w)|\le C{|a^{k,1}|\over \tau_{k, 2}}\le {C\over \tau_{k, 1}\tau_{k,2}}\leqno(3.39)
$$
and
$$
-T_2 T_3 A^2_{k,1} [{\d f_2\over \d \zbar_3}]
=T_2 A^2_{k, 1} B_{2,3} [f_2]+ T_3 A^2_{k,1}B_{3, 2}[f_3]+T_3 A^3_{2, k ,1} [f_3]\leqno(3.40)
$$
By (3.22), one has
\begin{eqnarray*}
\lefteqn{(4.41) \quad \tau_{2,3} \tau_{k, 1} \tau_{k,2}}\\
&\qquad \quad\ge & |w_1-z_1|^{2-\epsilon} |w_k-z_k|^\epsilon |w_k-z_k|^{2-2\epsilon}|w_2-z_2|^{2\epsilon} 
|w_2-z_2|^{2-3\epsilon} |w_3-z_3|^{3\epsilon}\\
&\qquad\quad =& |w_3-z_3|^{3\epsilon} \ell_{1, 2, k}(\epsilon).
\end{eqnarray*}
Applying  the inequality (4.41) and estimate (3.39), one has
$$
 \Big|{k(z_3, w_3)\over \tau_{2,3}} a^{2, k, 1}\Big |
\le {C  \over |w_3- z_3|  \tau_{2, 3} \tau_{k, 1} \tau_{k, 2}}
\le
{C\over |w_3-z_3|^{1+3\epsilon} \ell_{1,2, k}(\epsilon) },\leqno(3.42)
$$
$$
|k(z_3, w_3) {a^{k, 1} \over \tau_{2, k}}b^{3,2}|
\le {C\over |w_3-z_3| \tau_{2,3} \tau_{k,1} \tau_{2, k}}\le {C\over |w_3-z_3|^{1+3\epsilon} \ell_{1, 2, k}(\epsilon)}\leqno(3.43)
$$
and, similarly
$$
|k(z_2, w_2) {a^{k, 1} \over \tau_{2,k} }b^{2,3}|
\le {C\over |w_2-z_2| \tau_{2,3} \tau_{k,1} \tau_{2, k}}\le {C\over |w_2-z_2|^{1+3\epsilon} \ell_{1, 3, k}(\epsilon)} .\leqno(3.44)
$$
Therefore, combining the above estimates, the integral kernel of integral operators (3.40) can be written as
$T^\ell _{i, j, k}[f_\ell]  $ with integral kernel $E^\ell_{i,j,k}$  for any distinct $i, j, k, \ell\in \{1,2,\cdots, n\}$. Moreover,  $E^\ell_{i,j,k}$ satisfies the estimate
$$
|E^\ell _{i, j, k}|\le {C\over |w_\ell-z_\ell|^{1+3\epsilon} \ell_{i,j, k}(\epsilon)}.\leqno(3.45)
$$
Therefore, Theorem 3.2 is proved when $n=4$, it follows similarly when $n>4$ from all cases have been discussed above.
\epf
\smallskip

For any $n\in \NN$, we define: $\NN_n=\{1,2,\cdots, n\}$. 

\begin{proposition} For any $k\in \NN_n$ and $I=\{i_1,\cdots, i_m\} \subset \NN_n\setminus \{k\}$. Then
$T^k_I: L^p(\Omega^n)\to L^p(\Omega^n)$ is bounded and 
$$
\|T^k_I\|_{L^p(\Omega^n)\to L^p(\Omega^n)}\le C \|f  \|_{L^p(\Omega^n)},
\quad
\hbox{for all }\ 1\le  p\le \infty.
$$
\end{proposition}
\proof Since $T^k_I[g]=\int_{\Omega^\ell } E^k_I(z, w) g(w) dA(w_k, w^I)$ with $I=(i_1, \cdots, i_m)$
$$
|E^k_I(z,w)| \le {C\over |w_k-z_k|^{1+m\epsilon} \ell_I (\epsilon)}.
$$
Then
$$
\int_{\Omega^n} |E^k_I(z,w)| dv(w) \le {C\over \epsilon^n}  \ \hbox{ and }\  \int_{\Omega^n} |E^k_I(z,w)| dv(z) \le {C\over \epsilon^n}.
$$
By the Schur's lemma, one has 
$$
\|T^k_{I}\|_{L^p\to L^p} \le {C\over \epsilon^n},\quad 1<p<\infty.
$$
Since the constant $C \epsilon^{-n}$ is independent of $p$, by letting $p\to 1^+$ and then $p\to + \infty$, we have proved  the proof of the proposition.\epf
\medskip

As a corollary of Theorem 3.2 and Proposition 3.3, one has 
\begin{theorem} Let $f=\sum_{j=1}^n f_j d\zbar_j \in C^1_{(0,1)}(\Omegabar^n)$ be $\dbar$-closed. For $1\le j\le n$, there is a scalar constant $C$ such that
$$
\|T_j P_{j-1} \cdots P_1 P_0 f_j\|_{L^p(\Omega^n) }\le C \sum_{k=1}^j  \|f_k \|_{L^p(\Omega^n)}, \leqno(3. 46)
$$
for any $1\le p\le \infty$.
\end{theorem}

\section{Proof of Theorem 4.1} 

\subsection{Approximation}

\begin{theorem} Let $\Omega$ be a bounded simply connected domain in $\CC$ with $C^{1,\alpha}$ boundary for some $\alpha>0$. For any $1<p<\infty$ anf
$f\in L^p_{(0,1)}(\Omega^n)$ be $\dbar$-closed, then there is a  $\dbar$-closed squence $\{f_m\}_{m=1}^\infty \subset C^1_{(0,1)}(\Omegabar^n)$
such that 
$$
\lim_{m \to\infty} \|f_m-f\|_{L^p_{(0,1)}}=0.\leqno(4.1)
$$
\end{theorem}

\proof When $\Omega$ is the unit disk $D$, let $\chi^j\in C^\infty_0(D)$ be nonnegative and $\int_{D}\chi^j dA=1$. 
Let $\chi^j_\epsilon=\chi^j(z/\epsilon) \epsilon^{-2}$ and $\chi_\epsilon(z)=\chi^1_\epsilon \cdots \chi^n_\epsilon$
on $D^n$. The proof for this case  is very simple. For any $0<r<1$ and
$\epsilon=(1-r)/2$, since $f_r(z)=f(rz)$ is $\dbar$-closed in $D(0, 1/r)$ and then
$$
F_r(z)=f_r* \chi_\epsilon \in C^\infty_{(0,1)}(\overline{D}^n) \leqno(4.2)
$$
is $\dbar$-closed in $D^n$ and
$$
\|F_r-f\|_{L^p_{(0,1)}(D^n)}\to 0 \leqno(4.3)
$$
as $r\to 1^-$ and any $p\in (1,\infty)$. This argument remains true when $\Omega$ is a simply connected domain
in $\CC$ with $C^{1,\alpha}$ boundary for any $0<\alpha<1$. Let $\phi: \Omega\to D$ be a biholomorphic mapping. Then
$\phi\in C^{1,\alpha}(\Omegabar)$, and $\Omega=\phi^{-1}(D)$, with slightly modification of the unit disc case, one can
similarly prove the theorem.
 \epf

\medskip

Now we are ready to prove Theorem 1.1 when $\Omega$ is bounded simply connected with $C^{1,\alpha}$ boundary. 

\subsection{Proof of Theorem 1.1 when $\Omega$ is simply connected}

\proof For any $1<p<\infty$, by Theorem 4.1, there is a sequence $\{f_m \}_{m=1}^\infty\subset C^1_{(0,1)}(\Omegabar)$ which are
$\dbar$-closed such that
$$
\lim_{m\to\infty} \|f_m-f\|_{L^p_{(0,1)}(\Omega)}=0.\leqno(4.4)
$$
By estimations obtained in Section 3, one has that
$$
\dbar S[f_m]=f_m \leqno(4.5)
$$
and $S[f_m]$ is a canonical solution. Moreover,
$$
\lim_{m\to\infty} \|S[f_m]-S[f]\|_{L^p(\Omega^n)}=0. \leqno(4.6)
$$
For $2< p<\infty$, by Theorem 2.5, one has
\begin{eqnarray*}
\lefteqn{\|S[f]\|_{L^p(\Omega^n)}}\\
&\le & \|S[f_m]\|_{L^p (\Omega^n) } + \|S[f_m]-S[f]\|_{L^p (\Omega^n)}\\
&\le& C \|f_m\|_{L^p_{(0,1)}(\Omega^n)} + \|S[f_m]-S[f]\|_{L^p (\Omega^n)}\\
&\le & C\|f\|_{L^p_{(0,1)}(\Omega^n)} +C \|f_m-f\|_{L^p_{(0,1)}(\Omega^n)} + \|S[f_m]-S[f]\|_{L^p (\Omega^n)},
\end{eqnarray*}
where $C$ is a constant depends neither on  $m$ nor $p$.  Let $m\to \infty$, one has
$$
\|S[f]\|_{L^p_{(0,1)}(\Omega^n)}\le C\|f\|_{L^p_{(0,1)}(\Omega^n)}, \quad 2< p<\infty. \leqno(4.7)
$$
Letting $p\to +\infty$, one has
$$
\|S[f]\|_{L^\infty_{(0,1)}(\Omega^n)}\le C\|f\|_{L^\infty_{(0,1)}(\Omega^n)}.\leqno(4.8)
$$
The proof of Theorem 1.1 is complete when $\Omega$ is simply connected with $C^{1,\alpha}$ boundary.

\subsection{Proof of Theorem 1.1 for $\Omega$ satisfying the UEBC}

Since $\Omega$ is a bounded domain in $\CC$ with piecewise $C^1$ boundary and satisfies a uniform exterior ball condition
(of radius $r$), there is a sequence of domains $\Omega_\ell$ with piecewise $C^1$  boundary and
  satisfying the same
 uniform ball condition (of radius $r/2$) for all $\ell \ge 1$. Moreover,  
$$
\Omega_\ell\subset \Omegabar_\ell \subset \Omega_{\ell+1}\subset  \Omegabar_{\ell+1} \subset \Omega 
\quad\hbox{and}\quad
\lim_{\ell \to \infty} \Omega_{\ell} =\Omega. \leqno(4.9)
$$
Note, here we choose $\Omega_\ell$ so that the constant in Theorem 2.1 on the Green's function estimates on $\Omega_\ell$ is uniformly for all $\ell\ge 1$.

Notice that
$$
f*\chi_\epsilon \in C^\infty_{(0,1)}(\Omega_\ell^n)\leqno(4.10)
$$
is $\dbar$-closed in $\Omega_\ell$ if $\epsilon<\dist(\d \Omega_\ell, \d \Omega)$/n.  By the argument in Section 4.2, we have
$$
\|S_\ell[f] \|_{ L^p(\Omega_\ell^n)} \le C\|f\|_{L^p_{(0,1)} (\Omega_\ell ^n)},\quad \hbox{for } 2< p\le \infty, \leqno(4.11)
$$
where $C$ is a constant depend neither on $p$ nor  $\ell$. For any $1<p<\infty$, since the unit ball is weakly compact  in $L^p(\Omega_\ell)$, there is a subsequence
$\{S_{\ell_j}[f]\}_{j=1}^\infty $ converges to a function in $L^p(\Omega)$, denoted by $\tilde{S}[f]$ weakly on $L^p(\Omega_\ell)$ for any $\ell\ge 1$. Thus,
$$
\|\tilde{S}[f]\|_{L^p(\Omega_\ell^n)} \le C\|f\|_{L^p_{(0,1)}(\Omega_\ell^n)}\le C\|f\|_{L^p_{(0,1)}(\Omega^n)},\quad \ell\ge 1.\leqno(4.12)
$$
This implies that $\tilde{S}[f]\in L^p(\Omega^n)$ and
$$
\|\tilde{S}[f]\|_{L^p(\Omega^n)} \le C\|f\|_{L^p_{(0,1)}(\Omega^n)}. \leqno(4.13)
$$
By the uniqueness of weak limit for each $L^p(\Omega^n)$, one has $S[f]=\tilde{S}[f]$  for all $p\in (2,\infty)$. Since $C$ in (4.13) does not depend
 on $p$, letting $p\to\infty$, one has
$$
\|\tilde{S}[f]\|_{L^\infty(\Omega)^n} \le C\|f\|_{L^\infty_{(0,1)} (\Omega^n)}.\leqno(4.14)
$$
Since $S_\ell [f]$ is the canonical solution for $\dbar u=f$ in $\Omega_\ell$, it is easy to check $\dbar \tilde{S}[f]=f$ in $\Omega$ in the sense of distribution.
Moreover, for any $h\in L^2(\Omega)$, one has
$$
\int_{\Omega^n} \tilde{S}[f] \hbar(z) dv(z)=\lim_{\ell\to \infty} \int_{\Omega_\ell^n} S_\ell [f] \hbar(z) d(z)=0.\leqno(4.15)
$$
Therefore, $\tilde{S}[f]$ is the canonical solution of $\dbar u=f$ in $\Omega$. So, $S[f]=\tilde{S}[f],$  the proof is complete
when $\Omega$ satisfies a uniform ball condition.  Therefore, combining Sections 4.2 and 4.3, the proof of Theorem 1.1 is complete.\epf

\section{Remarks}

For any $\alpha \in [0, 1)$, we choose $\epsilon$ such that $(n+1) \epsilon=1-\alpha$. Thus, 
by the definition of $E^\ell _I$, one has  $|I|\le n-1$ and
$$
d_\Omega (w_k)^{-\alpha} |E^k_I(z,w)|
\le  {C\over |w_k-z_k|^{1+(n-1)\epsilon} d_\Omega (w_k)^{1-n\epsilon} \ell_I (\epsilon)}.\leqno(5.1)
$$
Therefore,  if  $ 1<p'\le {4-\epsilon\over 4-2\epsilon}$, then 
$$
\int_{\Omega^{\ell +1}}  \Big(d_\Omega (w_k)^{-\alpha} |E^k_I(z,w) |\Big)^{p'}  dA(w_k) dv(w_I)\le  {C \over \epsilon^{n}}. \leqno(5.2)
$$
This implies that
$$
|\int_{\Omega^{\ell +1}}  d_\Omega (w_k)^{-\alpha} E^k_I(z,w) f_k(w)  dA(w_k)dv(w_I) |\le  ({C \over \epsilon^{n}})^{1/p'} \|f_k\|_{L^p(\Omega^{\ell+1})}
$$
for all $p\ge {4-\epsilon \over \epsilon}$. Therefore,
$$
(5.3)\quad \Big\|\int_{\Omega^{\ell +1}}  d_\Omega (w_k)^{-\alpha} E^k_I(z,w) f_k(w)  dA(w_k)dv(w_I ) \Big \|_{L^p(\Omega^n )}
\le {C \over \epsilon^{n}} \|f_k\|_{L^p(\Omega^n)},
$$ 
for  all $p\ge {4-\epsilon \over \epsilon}$. Therefore, by (5.3) and arguments given in Section 4, we have proved the following theorem.

\begin{theorem} Let $\Omega$ be an admissible domain in $\CC$ and let $f=\sum_{j=1}^n f_j d\zbar_j \in L^\infty_{(0,1)}(\Omega^n)$ be $\dbar$-closed. Then
there is a scalar constant $C$ such that
$$
\|S[f]\|_{L^\infty (\Omega^n) }\le {C \over (1-\alpha)^n }\sum_{k=1}^n  \| d_\Omega (z_k)^\alpha f_k(z) \|_{L^\infty (\Omega^n)},  \leqno(5.4)
$$
for any $0<\alpha<1$.
\end{theorem}

\fontsize{11}{11}\selectfont
 \bigskip
 
\noindent Department of Mathematics, University of California, Irvine, CA 92697-3875, USA

 \bigskip
\noindent Email addresses:  \ sli@math.uci.edu

\end{document}